\documentclass{amsart}

\newtheorem{theorem}{Theorem}
\newcommand{\bt}{\begin{theorem}}
\newcommand{\et}{\end{theorem}}

\newtheorem{lemma}{Lemma}
\newcommand{\bl}{\begin{lemma}}
\newcommand{\el}{\end{lemma}}

\newtheorem{corollary}{Corollary}
\newcommand{\bc}{\begin{corollary}}
\newcommand{\ec}{\end{corollary}}

\DeclareMathOperator{\qqand}{\qquad\text{and}\qquad}

\newcommand{\N}{\ensuremath{ \mathbf N }}

\newcommand{\mce}{\ensuremath{ \mathcal E}}
\newcommand{\mcf}{\ensuremath{ \mathcal F}}
\newcommand{\mcg}{\ensuremath{ \mathcal G}}
\newcommand{\mcm}{\ensuremath{ \mathcal M}}
\newcommand{\mcw}{\ensuremath{ \mathcal W}}

\newcommand{\beq}{\begin{equation}}
\newcommand{\eeq}{\end{equation}}

\usepackage{amsmath,amssymb,amsthm}

\title{A new class of minimal asymptotic bases}
\author{Melvyn B. Nathanson}
\address{Lehman College (CUNY), Bronx, NY 10468}
\email{melvyn.nathanson@lehman.cuny.edu}
\date{\today}



\thanks{Supported in part by a grant from the PSC-CUNY Research Awards Program.}

\begin{document}

\begin{abstract}
A set $A$ of nonnegative integers is an asymptotic basis of order $h$ if every sufficiently 
large integer can be represented as the sum of $h$ not necessarily distinct elements of $A$.  
An asymptotic basis $A$ is minimal if 
removing any element of $A$ destroys every representation of infinitely many integers.  
In this paper,  a new class of minimal asymptotic bases is constructed.  
\end{abstract}

\maketitle

\section{\mcg-adic asymptotic bases}

Let $\N_0 = \{0,1,2,3,\ldots\}$ be the set of nonnegative integers and 
let $h$ be a positive integer.  
Let $A, A_1,\ldots, A_h$ be subsets of $\N_0$.  
We define the \emph{sumset} 
\[
A_1 + \cdots + A_h = \{a_1 + \cdots + a_h: a_i \in A_i \text{ for all } i  = 1,2,\ldots, h \} 
\]
and the \emph{$h$-fold sumset} 
\[
hA = \underbrace{A + \cdots + A}_{\text{$h$ summands}} = \{a_1 + \cdots + a_h: 
a_i \in A \text{ for all } i =1,2,\ldots, h\}.
\]
The set $A$ is a \emph{basis of order $h$} if every nonnegative integer can be represented 
as the sum of $h$ not necessarily distinct elements of $A$, that is, if $hA = \N_0$.  
The set $A$ is an \emph{asymptotic basis of order $h$} if $hA$ contains all sufficiently large integers.   
An asymptotic basis of order $h$ is \emph{minimal} if no proper subset of $A$ 
is an asymptotic basis of order $h$.  
Thus, if $A$ is a minimal asymptotic basis of order $h$, then,
for all $a \in A$, there are infinitely many integers $n$ such that removing $a$ from $A$ destroys 
every representation of $n$ as a sum of $h$ elements of $A$.  

Minimal asymptotic bases are extremal objects in additive number theory, and are related to 
the conjecture of Erd\H os and Tur\' an~\cite{erdo-tura41} that the representation function 
of an asymptotic basis of order $h$ must be unbounded.   

At this time there are few explicit constructions of minimal asymptotic bases.  
Nathanson~\cite{nath1974-11,nath1988-66} used a 2-adic construction to produce 
the first examples of minimal asymptotic bases.    
This method was extended to $g$-adically defined sets by  
Chen~\cite{chen12}, 
Chen and Chen~\cite{chen-chen11}, Chen and Tang~\cite{chen-tang18}, 
Jia~\cite{jia96}, Jia and Nathanson~\cite{nath1989-68}, 
Lee~\cite{lee93}, Li and Li~\cite{li-li16}, 
Ling and Tang~\cite{ling-tang15,ling-tang18}, 
Sun~\cite{sun19,sun21}, and 
Sun and Tao~\cite{sun-tao19}.
This paper constructs a new class of minimal asymptotic bases.

An \emph{interval of integers of length $t$} is a set of $t$ consecutive integers.  
For $u,v \in \N_0$ with $u \leq v$, the set 
\[
[u,v] = \{x\in \N_0: u \leq x \leq v\} 
\]
is an interval of integers of length $v - u + 1$.  
Thus,  $[1,h] = \{1,2,\ldots, h\}$.

A \emph{\mcg-adic sequence} is a strictly increasing sequence of positive integers 
$\mcg = (g_i)_{i=0}^{\infty}$ 
such that $g_ 0 = 1$ and $g_{i-1}$ divides $g_i$ for all $ i \geq 1$.   
Let $(d_i)_{i=1}^{\infty}$ be the sequence  of  positive integers defined by 
\[
d_i = \frac{g_i}{g_{i-1}}.
\]
For all $i \geq 1$, we have 
\beq                                         \label{MinBasis:g1}
d_i \geq 2
\eeq
and
\beq                                         \label{MinBasis:g2}
g_i = d_1 d_2 \cdots d_i.  
\eeq
For all $i$ and $j$ with $0 \leq i <  j$, we have 
\beq                                        \label{MinBasis:g3}
\frac{g_{i+j}}{g_i} = d_{i+1}  d_{i+2} \cdots d_{i+j}.
\eeq
Every positive integer $n$ has a unique \textit{$\mcg$-adic representation}  
\[
n = \sum_{j=0}^{\infty} x_j g_j 
\]
where 
\[
x_j \in [0, d_{j+1}-1] 
\]
 for all $j \in \N_0$ and $x_j = 0$ for all sufficiently large $j$ 
(Nathanson~\cite{nath2014-150,nath2017-172}).  
Equivalently, for every positive integer $n$, there is a unique nonempty finite set $F \subseteq \N_0$ 
and a unique set $\{x_j : j  \in F \}$ such that 
\beq                             \label{MinBasis:g-adic}
n = \sum_{j \in F} x_j g_j 
\eeq
where $x_j \in  [1, d_{j+1}-1]$ for all $j \in F$.

For every integer $g \geq 2$, the usual $g$-adic representation uses the 
$\mcg$-adic sequence $\mcg = \left(g^i\right)_{i=0}^{\infty}$ with 
quotients $d_i = g$ for all $i \geq 1$.

\bl                                             \label{MinBasis:lemma:Gbound} 
Let $\mcg = (g_i)_{i=0}^{\infty}$ be a \mcg-adic system.  
If $n =  \sum_{j \in F} x_j g_j $ is a positive integer 
with $x_j \in  [1, d_{j+1}-1]$ for all $j \in F$, 
then 
\[
g_M \leq n < g_{M+1}
\]
if and only if 
\[
\max(F) = M.  
\]
\el

\begin{proof}
If $\max(F) = M$, then $F \subseteq [0,M]$ and 
\begin{align*}
g_M & \leq x_M g_M \leq n \leq  \sum_{j \in F} (d_{j+1}-1) g_j \\
&  \leq  \sum_{j =0}^M (d_{j+1}-1) g_j =  \sum_{j =0}^M (g_{j+1}- g_j) \\
& = g_{M+1}-1 < g_{M+1}.  
\end{align*}

Conversely, if $g_M \leq n < g_{M+1}$ and $\max(F) = M'$, 
then the inequality 
\[
g_{M'} \leq x_{M'} g_{M'} \leq n < g_{M+1}
\]
implies $M' \leq M$.  If $M' \leq M-1$, then 
\[
n \leq \sum_{j =0}^{M'} (d_{j+1}-1) g_j < g_{M'+1} \leq g_M
\]
which is absurd.  Therefore, $g_M \leq n < g_{M+1}$ implies 
$\max(F) = M$. 
This completes the proof. 
\end{proof}

Let $W$ be a nonempty set of nonnegative integers, 
and let $\mcf^*(W)$ be the set of all nonempty finite subsets of $W$.  
Let $\mcg = (g_i)_{i=0}^{\infty}$ be a \mcg-adic sequence. 
We define the set of positive integers 
\[
A_{\mcg}(W) = \left\{ \sum_{j\in F} x_j g_j : F \in \mcf^*(W) \text{ and }  x_j \in [1,d_{j+1} -1] \right\}.
\]
Note that $0 \notin A_{\mcg}(W)$ because $\emptyset \notin  \mcf^*(W)$.

Let $h \geq 2$.  
A \emph{partition} of $\N_0$ is a sequence $\mcw = (W_i)_{i=0}^{h-1}$ 
of nonempty pairwise disjoint sets such that 
\[
\N_0 = W_0 \cup W_1 \cup \cdots \cup W_{h-1}.
\]

\bt                  \label{MinBasis:theorem:basis}
Let $h \geq 2$ and let $\mcw = (W_i)_{i=0}^{h-1}$ be a partition of $\N_0$.  
Let $\mcg = (g_i)_{i=0}^{\infty}$ be a \mcg-adic sequence. 
The set 
\[
A_{\mcg}(\mcw) = \bigcup_{i=0}^{h-1} A_{\mcg}(W_i)
\]
is an asymptotic basis of order $h$ with $h$-fold sumset  
\[
hA_{\mcg}(\mcw) = \{n \in \N_0: n \geq h \}.
\]
\et

\begin{proof} 
The smallest integer in the set  $A_{\mcg}(\mcw)$ is $1 = 1 \cdot g_0$.    
It follows that $h \in hA_{\mcg}(\mcw)$ but 
$[0,h-1] \cap hA_{\mcg}(\mcw) = \emptyset$.  

Every positive integer $n$ has a unique \mcg-adic representation $n = \sum_{j \in F} x_j g_j$, 
where $F$ is a nonempty finite set of nonnegative integers and $x_j \in [0, d_{j+1}-1] $.  
For $i \in [0,h-1]$, let 
\[
F_i = F  \cap W_i   
\]
and 
\[
n_i = \sum_{j\in F_i} x_j g_j. 
\]
If   $F_i = \emptyset$, then $n_i = 0 \notin  A_{\mcg}(\mcw)$.  If $F_i \neq \emptyset$, then 
\[
n_i = \sum_{j\in F_i} x_j g_j \in A_{\mcg}(W_i) \subseteq A_{\mcg}(\mcw).   
\]
Let 
\[
L = \left\{ i\in [0,h-1] :  F_i \neq \emptyset   \right\} = \left\{ i\in [0,h-1] :  n_i \geq 1  \right\} 
\] 
and 
\[
|L| = \ell_0. 
\]
We have  $\ell_0 \in [1,h]$ and 
\[
n = \sum_{i\in L} n_i \in \ell_0  A_{\mcg}(\mcw). 
\]
Let $\ell$ be the largest integer such that  $\ell \leq h$ and 
$n \in \ell  A_{\mcg}(\mcw)$.   We must prove that $\ell = h$. 

If $n \in \ell  A_{\mcg}(\mcw)$, then there exist integers 
$n_1,\ldots, n_{\ell} \in A_{\mcg}(\mcw)$ 
such that 
\[
n = n_1 + \cdots  + n_{k-1} +  n_k +  n_{k+1} + \cdots + n_{\ell}. 
\]
For each $i \in [1,\ell]$ there is an integer $s_i \in [0,h-1]$ and  a set $F_{s_i} \in \mcf^*(W_{s_i})$ 
such that $n_i$ has the \mcg-adic representation 
\[
n_i = \sum_{j \in F_{s_i}} x_j g_j \in A_{\mcg}(W_{s_i}). 
\] 

Suppose that $\ell < h$. If $|F_{s_k}| \geq 2$ for some $k \in [1, \ell]$, 
then there are nonempty sets $F'_{s_i}$ and $F''_{s_i}$ such that 
\[
F_{s_k} = F'_{s_k} \cup F''_{s_k}  \qqand F'_{s_k} \cap F''_{s_k} = \emptyset.
\]
The  integers 
\[
n'_k = \sum_{j \in F'_{s_k}} x_j g_j  \in A_{\mcg}(W_{s_k}) 
 \qqand n''_k = \sum_{j \in F''_{s_k}} x_j g_j  \in A_{\mcg}(W_{s_k}) 
\]
satisfy 
\[
n_k = n'_k + n''_k 
\]
and so 
\[
n = n_1 + \cdots n_{k-1} +  n'_k + n''_k +  n_{k+1} + \cdots + n_{\ell} \in (\ell+1) A_{\mcg}(\mcw). 
\]
This contradicts the maximality of $\ell$, and so $|F_{s_i}| = 1$ for all $i \in [1,\ell]$ and
\[
n_i = x_{j_i} g_{j_i}
\]
for some $g_{j_i} \in W_{s_i}$ and $x_{j_i} \in [0, d_{j_i+1}-1]$.  

If $x_{j_k} \geq 2$ for some $k \in [1,\ell]$, then 
\[
g_{j_k} \in  A_{\mcg}(W_{s_k}) \qqand  (x_{j_k} -1) g_{j_k} \in  A_{\mcg}(W_{s_k})  
\]
and 
\[
n_k = g_{j_k} + (x_{j_k} -1) g_{j_k}.
\]
It follows that 
\[
n = n_1 + \cdots n_{k-1} +  g_{j_k} + (x_{j_k} -1) g_{j_k} +  n_{k+1} + \cdots + n_{\ell} \in (\ell+1) A_{\mcg}(\mcw), 
\]
which again contradicts the maximality of $\ell$.  Therefore, $x_{j_i} = 1$ for all $i \in [1,\ell]$, and 
\[
n = g_{j_1}  + \cdots + g_{j_k}   + \cdots + g_{j_{\ell}}.    
\]

If $j_k \geq 1$ for some $k \in [1,\ell]$, then 
\[
g_{j_k} = g_{j_k-1} + (d_{j_k}-1) g_{j_k-1} 
\]
and 
\[
n = g_{j_1}  + \cdots + g_{j_k-1} + (d_{j_k}-1) g_{j_k-1} + \cdots + g_{j_{\ell}}     
\in (\ell+1) A_{\mcg}(\mcw), 
\]
which also contradicts the maximality of $\ell$.  
Therefore, $j_k = 0$  and $g_{j_k} = g_0 = 1$ for all $k \in [1,\ell]$, and so 
\[
n = \underbrace{g_0  + \cdots + g_0}_{\text{$\ell$ summands}}   
= \underbrace{1+ \cdots + 1 + \cdots + 1}_{\text{$\ell$ summands}} 
= \ell \leq h-1, 
\]
which is absurd.  
This completes the proof. 
\end{proof} 

\bt                         \label{MinBasis:theorem:Abasis}
Let $h \geq 2$ and let $\mcw = (W_i)_{i=0}^{h-1}$ be a partition of $\N_0$.  
The set 
\[
A = \{0\} \cup A_{\mcg}(\mcw) 
\]
is a basis of order $h$ but not a minimal asymptotic basis of order $h$. 
\et

\begin{proof}
Because $0 \in A$ and $1 \in A_{\mcg}(\mcw) \subseteq A$, we have   
\[
\ell = (h-\ell)\cdot 0 + \ell \cdot 1 \in hA
\]
 for all $\ell \in [0,h-1]$. 
By Theorem~\ref{MinBasis:theorem:basis}, 
\[
 \{n \in \N_0: n \geq h \} = hA_{\mcg}(\mcw) \subseteq hA
 \]
 and so $A$ is a basis of order $h$.  

The set $A$ is not a minimal asymptotic basis of order $h$ 
because $0 \in A$ and the removal of 0 from $A$ gives the set $A \setminus \{0 \} = A_{\mcg}(\mcw)$,
which is still an asymptotic basis of order $h$. 
This completes the proof.  
\end{proof}

\section{Minimal asymptotic bases}

The following lemma generalizes a result of Jia~\cite{jia96}.

\bl                   \label{MinBasis:lemma:ineqM} 
Let $\mcg = (g_i)_{i=0}^{\infty}$ be a \mcg-adic sequence.
Let $(u_i)_{i=1}^p$ be a strictly increasing finite sequence of nonnegative integers, 
and let $(v_j)_{j=1}^q$ be a finite sequence of not necessarily distinct nonnegative integers.  
Let 
\[
x_i \in [1, d_{u_i+1}-1]
\]
 for all $i \in [1,p]$ and 
 \[
 y_j \in [1, d_{v_j+1}-1] 
 \]
 for  all $j \in [1,q]$.   
If 
\beq                                 \label{MinBasis:ineqM-n}
n = \sum_{i=1}^p x_i g_{u_i} = \sum_{j=1}^q y_j g_{v_j} 
\eeq
then 
\beq                                 \label{MinBasis:ineqM}
\sum_{u_i\leq u_k} x_i g_{u_i}   \leq  \sum_{v_j \leq u_k} y_j g_{v_j}  
\eeq
for all $k \in [1, p]$. 
\el

\begin{proof}
Because the sequence  $(u_i)_{i=1}^p$is  strictly increasing, 
Lemma~\ref{MinBasis:lemma:Gbound}  implies 
\begin{align*} 
\sum_{u_i\leq u_k} x_i g_{u_i}  =  \sum_{i=1}^k x_i g_{u_i}  < g_{u_{k+1}} 
\end{align*} 
for all $k \in [1,p]$.    
Choosing $k = p$ gives 
\[
n  = \sum_{i=1}^p x_i g_{u_i} < g_{u_{p+1}}.
\]
Relation~\eqref{MinBasis:ineqM-n} implies 
\[
g_{v_j} \leq n < g_{u_p + 1} 
\]
an so  $v_j \leq u_p$ for all $j \in [1,q]$.  
This implies~\eqref{MinBasis:ineqM} for $k=p$.  

In the sequence $\mcg = (g_i)_{i=0}^{\infty}$, the integer  $g_i$ divides $g_j$ 
for all $i \leq j$.  
Let $k \in [1,p-1]$.  
Because the  sequence $(u_i)_{i=1}^p$ is strictly increasing, for all $i \in [k+1,p]$ we have 
\[
u_k < u_k+1 \leq u_{k+1} \leq u_i 
\]
and $g_{u_k+1}$ divides $g_{u_i}$, that is, 
\[
g_{u_i} \equiv 0 \pmod{g_{u_k+1}}.
\]
If $v_j \geq u_k +1$, then 
\[
g_{v_j} \equiv 0 \pmod{g_{u_k+1}}.
\]
Rearranging~\eqref{MinBasis:ineqM-n}, we obtain 
\begin{align}                                   \label{MinBasis:ineqM-cong}  
\sum_{u_i\leq u_k} x_i g_{u_i} -  \sum_{v_j \leq u_k} y_j g_{v_j} 
&  =   \sum_{v_j \geq u_k +1} y_j g_{v_j} - \sum_{u_i > u_k} x_i g_{u_i}   \\
&  \equiv 0 \pmod{g_{u_k+1}}.    \nonumber 
 \end{align}    
If  
\[
\sum_{u_i\leq u_k} x_i g_{u_i}  >  \sum_{v_j \leq u_k} y_j g_{v_j} 
\]
then Lemma~\ref{MinBasis:lemma:Gbound} gives 
\[
0 <   \sum_{i=1}^k x_i g_{u_i}  - \sum_{v_j \leq u_k} y_j g_{v_j} 
\leq  \sum_{i=1}^k x_i g_{u_i} < g_{u_{k+1}}. 
\]
This inequality contradicts congruence~\eqref{MinBasis:ineqM-cong}. 
This completes the proof.    
\end{proof}

\bt                              \label{MinBasis:theorem:NewMinBasis}
Let $h \geq 2$ and let  $t$ be an integer such that 
\beq             \label{MinBase:t}
t \geq 1 +  \frac{\log h}{\log 2}.
\eeq 
Let $\mcw = (W_i)_{i=0}^{h-1}$ be a partition of $\N_0$  such that, 
for all $i \in [0 ,h-1]$, there is an infinite set $\mcm_i$ of positive integers  such that 
\beq             \label{MinBase:t-interval}
[M_i -t+ 1 , M_i] \subseteq W_i
\eeq
for all $M_i \in \mcm_i$.   
Let $\mcg = (g_i)_{i=0}^{\infty}$ be a \mcg-adic sequence.  
The set 
\[
A_{\mcg}(\mcw) = \bigcup_{i=0}^{h-1} A_{\mcg}(W_i)
\]
is a minimal asymptotic basis of order $h$.  
\et

\begin{proof}
By Theorem~\ref{MinBasis:theorem:basis}, the set $ A_{\mcg}(\mcw) $ 
is an asymptotic basis of order $h$. 

Let $a \in A_{\mcg}(\mcw)$.  
Without loss of generality, we can assume that $a = a_0 \in A_{\mcg}(W_0)$ 
and 
\beq             \label{MinBase:a0}
a_0 = \sum_{j \in F_0} x_{0,j} g_j
\eeq
where $F_0 \in \mcf^*(W_0)$ and $x_{0,j} \in [1, d_{j+1}-1]$ for all $j \in F_0$.
Let $M_0 = \max(F_0)$.  By Lemma~\ref{MinBasis:lemma:Gbound},  
\[
g_{M_0} \leq a_0 < g_{M_0+1}. 
\]  

For all $i \in [1,h-1]$, choose an integer $M_i$ in the infinite set $\mcm_i$  
such that  
\beq             \label{MinBase:MiMot}
M_i \geq M_0 + t 
\eeq
and let  
\beq             \label{MinBase:ai}
a_i = \sum_{\substack{j \in W_i \\ j < M_0 }} (d_{j+1} - 1) g_j + g_{M_i} \in A_{\mcg}(W_i).
\eeq 
This is the \mcg-adic representation of $a_i$.  
Let 
\beq             \label{MinBase:nn}
n = \sum_{i=0}^{h-1} a_i 
= a_0 + \sum_{i=1}^{h-1} \sum_{\substack{j \in W_i \\ j < M_0 }} (d_{j+1} - 1) g_j 
+ \sum_{i=1}^{h-1} g_{M_i}.
\eeq
This is the \mcg-adic representation of $n$.  

Let 
\[
n = \sum_{i=0}^{h-1} b_{k_i} 
\]
be any representation of $n$ as the sum of $h$ elements of $A_{\mcg}(\mcw)$, 
where $k_i \in  [0,h-1]$ for all $i\in  [0,h-1]$ 
and $b_{k_i} \in A_{\mcg}(W_{k_i})$.   
We must prove that $b_{k_i} = a_0$ for some $i \in [0,h-1]$. 

Each integer $b_{k_i}$  is of the form 
\[
b_{k_i} = \sum_{j \in E_{k_i}} y_{i,j} g_j \in  A_{\mcg}(W_{k_i})   
\]
where $E_{k_i} \in \mcf^*(W_{k_i})$  
and $y_{i,j} \in [1,g_{j+1} - 1]$ for all $j \in E_{k_i}$. 
The uniqueness of the \mcg-adic representation implies that if $\{k_0,k_1,\ldots, k_{h-1}\} = [0,h-1]$,
then, after rearrangement, $k_i = i$ and $a_i = b_i$ for all $i \in [0,h-1]$.  

If $\{k_0, k_1,\ldots, k_{h-1}\} \neq [0,h-1]$, then there exists $s \in [0,h-1]$ 
such that $s \notin \{k_0,k_1,\ldots, k_{h-1}\}$.   
Suppose that $s \neq 0$. 
Recall that  
\[
M_s \geq M_0 + t 
\]
and 
\[
[M_s -t+1,M_s] \subseteq W_s. 
\]
 Because $k_i \neq s$ for all $i \in [0,h-1]$, we have  
 \[
[M_s - t + 1, M_s] \cap E_{k_i} 
\subseteq    W_s\cap W_{k_i}  = \emptyset.     
\] 
We construct the partition   
\[
E_{k_i} = E'_{k_i} \cup E''_{k_i}
\]
with 
\[
E_{k_i}' = \{j \in E_{k_i}: j \leq M_s - t\}  
\]
and 
\[
E_{k_i}'' = \{j \in E_{k_i}: j \geq M_s+1 \}.   
\] 
The sets $E_{k_i}' $ and $E_{k_i}''$ are not necessarily nonempty.  
Let 
\[
b_{k_i} = b'_{k_i} + b''_{k_i} 
\]
where 
\[
b'_{k_i} = \sum_{j \in E'_{k_i}} y_{i,j} g_j   \qqand   b''_{k_i} = \sum_{j \in E''_{k_i}} y_{i,j} g_j.  
\]
Note that $b'_{k_i} = 0$ if  $E'_{k_i} = \emptyset$ and $b''_{k_i} = 0$ if  $E''_{k_i} = \emptyset$.  

By Lemma~\ref{MinBasis:lemma:Gbound},  
\beq                              \label{MinBase:g5}
b'_{k_i} < g_{M_s -t +1} 
\eeq
and 
\[
b''_{k_i} = 0  \qquad\text{or}\qquad b''_{k_i} \geq g_{M_s+1}. 
\]
Let 
\beq                             \label{MinBase:g6}
n' = \sum_{i=0}^{h-1} b'_{k_i}  \qqand n'' = \sum_{i=0}^{h-1} b''_{k_i}.
\eeq
Recall inequalities~\eqref{MinBasis:g1}  and~\eqref{MinBase:t}: 
\[
d_i \geq 2 \qqand h \leq 2^{t-1}. 
\]
From~\eqref{MinBase:g6},~\eqref{MinBase:g5}, and~\eqref{MinBasis:g3},  we obtain 
\[
n' < hg_{M_s -t+1} \leq 2^{t-1} g_{M_s -t+1} \leq g_{M_s -t+1} \prod_{i=1}^{t-1}d_{M_s -t+i} =g_{M_s}. 
\]
Therefore, the \mcg-adic representation of $n'$ is of the form
\[
n' = \sum_{j =0}^{M_s -1} z_ j g_j 
\]
with $z_j \in [0,d_{j+1}-1]$.  Because 
\[
n'' = 0 \qquad\text{or}\qquad  n'' \geq g_{M_s+1}, 
\]
the \mcg-adic representation of $n''$ is of the form
\[
n'' = \sum_{j= M_s+1}^{\infty}  z_j g_j 
\]
with $z_j \in [0,d_{j+1}-1]$ and $z_j \geq 1$ for only finitely many $j$.  
Therefore, 
\[
n = n' + n'' = \sum_{i=0}^{ M_s -1 } z_i g_i + \sum_{i= M_s+1}^{\infty}  z_i g_i
\]
is the \mcg-adic representation of $n$.  
In this representation,  the coefficient of $g_{M_s}$ is 0, which contradicts the construction of $n$.
It follows that 
\[
[1,h-1] \subseteq \{k_0,k_1,\ldots, k_{h-1}\}.
\]
Renumbering the integers $b_i$, we can assume that $k_i = i$ 
and $b_i \in A_{\mcg}(W_i)$ for all $i \in [1,h-1]$.

We must prove that $k_0 = 0$, or, equivalently, that $b_0 \in A_{\mcg}(W_0)$.  
If not, then $b_0 \in A_{\mcg}(W_r)$ for some $r \in [1,h-1]$.
Because 
\[
M_0 = \max(F_0)  \in F_0 \subseteq W_0 
\]
we have  
\[
M_0 \notin \bigcup_{i=0}^{h-1} E_{k_i} \subseteq \bigcup_{i=1}^{h-1} W_i. 
\]    
Identity~\eqref{MinBase:a0} and Lemma~\ref{MinBasis:lemma:Gbound} give 
\[
g_{M_0} \leq a_0 = \sum_{j \in F_0} x_{0,j} g_j < g_{M_0 +1}. 
\]
From~\eqref{MinBase:ai} we have 
\begin{align*} 
n & = 
 a_0 + \sum_{i=1}^{h-1} \sum_{\substack{j \in W_i \\ j < M_0 }} (d_{j+1} - 1) g_j 
+ \sum_{i=1}^{h-1} g_{M_i} \\
& = \sum_{i=0}^{h-1} \sum_{ j \in E_i } y_{i,j} g_j. 
\end{align*} 
Summing only over terms $g_j$ with $j \leq M_0$ and 
applying Lemma~\ref{MinBasis:lemma:ineqM}, we obtain 
\begin{align*}
a_0 +  \sum_{i=1}^{h-1} \sum_{ \substack{j \in W_i \\ j < M_0 } } (d_{j+1}-1)  g_j 
& \leq \sum_{i=0}^{h-1} \sum_{\substack{j \in E_i \\ j \leq M_0}} y_{i,j} g_j 
=  \sum_{i=0}^{h-1} \sum_{\substack{j \in E_i \\ j < M_0}} y_{i,j} g_j \\
& =  \sum_{\substack{j \in E_0 \\ j < M_0}} y_{0,j} g_j 
 + \sum_{i=1}^{h-1} \sum_{\substack{j \in E_i \\ j < M_0}} y_{i,j} g_j \\
& < g_{M_0} +  \sum_{i=1}^{h-1}  \sum_{ \substack{j \in W_i \\ j < M_0 } } (d_{j+1}-1)  g_j \\ 
& \leq a_0  +  \sum_{i=1}^{h-1}  \sum_{ \substack{j \in W_i \\ j < M_0 } } (d_{j+1}-1)  g_j. 
\end{align*}
This is absurd, and so $b_0 \in A_{\mcg}(W_0)$. 
It follows that the integer $n$ defined by~\eqref{MinBase:nn} has a unique representation 
as the sum of $h$ elements of $A_{\mcg}(\mcw)$.  Therefore,  
\beq             \label{MinBase:NotIn}
n \notin h \left(A_{\mcg}(\mcw) \setminus \{a\} \right).
\eeq

For all $i \in [1,h-1]$, there are infinitely many integers $M_i \in \mcm_i$ 
with $M_i \geq M_0 + t$,
and so infinitely many positive integers $n$ satisfying~\eqref{MinBase:NotIn}.  
It follows that $A_{\mcg}(\mcw) \setminus \{a\}$ is not an asymptotic basis of order $h$ 
for all $a \in A_{\mcg}(\mcw)$.  
Equivalently, $A_{\mcg}(\mcw)$ is a minimal asymptotic basis of order $h$.   
This completes the proof.
\end{proof}

\bc                           
Let $\mcw = W_0 \cup W_1$ be a partition of $\N_0$  such that both $W_0$ and $W_1$ contain 
infinitely many pairs of consecutive integers.  
Let $\mcg = (g_i)_{i=0}^{\infty}$ be a \mcg-adic sequence.  
The set 
\[
A_{\mcg}(\mcw) =  A_{\mcg}(W_1) \cup A_{\mcg}(W_2)
\]
is a minimal asymptotic basis of order $2$.  
\ec 

\begin{proof}
This is the case $h=2$ of Theorem~\ref{MinBasis:theorem:NewMinBasis}.
\end{proof}

\section{Open problems} 
\begin{enumerate}
\item               \label{MinBasis:problem-a}
Let $h \geq 2$.  By Theorem~\ref{MinBasis:theorem:basis}, 
for every partition $\mcw = (W_i)_{i=0}^{h-1}$  of $\N_0$ and every \mcg-adic sequence, the set 
\[
A = \{0\} \cup A_{\mcg}(\mcw)
\]
is a basis of order $h$ that is not a minimal asymptotic basis of order $h$.
\begin{enumerate}
\item
Determine the set of all integers $a\in A$ such that $A \setminus \{a\}$ is an asymptotic basis of order $h$.   
\item
Determine the  set of all integers $a \in A$ such that $A\setminus \{a\}$ is a minimal asymptotic basis.  
\end{enumerate}

\item                  \label{MinBasis:problem-b}
Let $\mcg = (g_i)_{i=0}^{\infty}$ be a \mcg-adic sequence. Let $h \geq 2$.  
\begin{enumerate}
\item
Construct partitions $\mcw = (W_i)_{i=0}^{h-1}$ of $\N_0$ 
such that the set $A_{\mcg}(\mcw)$ is a minimal asymptotic basis of order $h$.  
\item
Construct partitions $\mcw = (W_i)_{i=0}^{h-1}$ of $\N_0$ 
such that the set $A_{\mcg}(\mcw)$ is not a minimal asymptotic basis of order $h$.  
\end{enumerate}

\item                     \label{MinBasis:problem-c} 
Let $(d_i)_{i=1}^{\infty}$ be a sequence of 2s and 3s.
Let  $\mcg = (g_i)_{i=0}^{\infty}$ be the \mcg-adic sequence defined by $g_0 =1$ 
and $g_i = \prod_{j=1}^i  d_i$ for $i \geq 1$.   
Consider problem~\eqref{MinBasis:problem-b} with respect to this \mcg-adic sequence.  
Of particular interest are the infinitely many \mcg-adic sequences $\mcg = (g_i)_{i=0}^{\infty}$ 
with quotients $\{ d_{2i-1}, d_{2i} \} = \{ 2,3 \}$ for all $i = 1,2,3,\ldots$.  In this case, 
$g_{2i} = 6^i $ for all $i$.  

\end{enumerate}

\def\cprime{$'$} \def\cprime{$'$} \def\cprime{$'$}
\providecommand{\bysame}{\leavevmode\hbox to3em{\hrulefill}\thinspace}
\providecommand{\MR}{\relax\ifhmode\unskip\space\fi MR }
\providecommand{\MRhref}[2]{%
  \href{http://www.ams.org/mathscinet-getitem?mr=#1}{#2}
}
\providecommand{\href}[2]{#2}

\end{document}